\documentclass[leqno]{article}
\usepackage{amssymb}
\usepackage{amsmath, enumerate, vmargin}

\makeatletter
\renewcommand\section{\@startsection {section}{1}{\z@}%
 {-3.5ex \@plus -1ex \@minus -.2ex}%
 {2.3ex \@plus.2ex}%
 {\center \normalfont\large\bfseries}}
\makeatother

\setmarginsrb{3cm}{3cm}{3cm}{3cm}{1mm}{6mm}{0mm}{10mm}

\newtheorem{thm}{Theorem}[section]

\newtheorem{cor}[thm]{Corollary}
\newtheorem{lem}[thm]{Lemma}
\newtheorem{defi}[thm]{Definition}
\newtheorem{remark}[thm]{Remark}
\newtheorem{example}[thm]{Example}
\newtheorem{pb}[thm]{Problem}

\newenvironment{rk}{\begin{remark}\rm}{\end{remark}}

\newcommand{\real}{{\mathbb R}}

\newcommand{\com}{{\mathbb C}}
\newcommand{\un}{1\mkern -4mu{\textrm l}}
\newcommand{\T}{{\mathbb T}}

\renewcommand{\L}{{\mathcal L}}
\renewcommand{\S}{{\mathcal S}}

\renewcommand{\a}{\alpha}
\renewcommand{\b}{\beta}

\newcommand{\Ga}{\Gamma}
\renewcommand{\d}{\delta}

\newcommand{\e}{\varepsilon}

\renewcommand{\l}{\lambda}

\newcommand{\f}{\varphi}
\renewcommand{\O}{{\Omega}}

\newcommand{\8}{\infty}
\newcommand{\el}{\ell}

\newcommand{\la}{\langle}
\newcommand{\ra}{\rangle}
\newcommand{\wt}{\widetilde}
\newcommand{\wh}{\widehat}
\newcommand{\n}{\noindent}
\newcommand{\pf}{\noindent{\it Proof.~~}}
\newcommand{\cqd}{\hfill$\Box$}
\newcommand{\be}{\begin{eqnarray*}}
\newcommand{\ee}{\end{eqnarray*}}
\newcommand{\beq}{\begin{equation}}
\newcommand{\eeq}{\end{equation}}

\numberwithin{equation}{section}

\begin{document}


\title{BMO functions and Carleson measures with values in uniformly
convex spaces}

\author{Caiheng Ouyang and Quanhua Xu}

\date{}

\maketitle


\begin{abstract}

 This paper studies the relationship between
vector-valued BMO functions  and the Carleson measures defined by
their gradients. Let $dA$ and $dm$ denote Lebesgue measures on the
unit disc $D$ and the unit circle $\T$, respectively. For $1<
q<\8$ and a Banach space $B$ we prove that there exists a positive
constant $c$ such that
  $$\sup_{z_0\in D}\int_{D}(1-|z|)^{q-1}\|\nabla f(z)\|^q
 P_{z_0}(z)\,dA(z)
 \le c^q\sup_{z_0\in D}\int_{\T}\|f(z)-f(z_0)\|^qP_{z_0}(z)\,dm(z)$$
holds for all trigonometric polynomials $f$ with coefficients in
$B$ iff $B$ admits an equivalent norm which is $q$-uniformly
convex, where
 $$P_{z_0}(z)=\frac{1-|z_0|^2}{|1-\bar{z_0}z|^2}\,.$$
The validity of the converse inequality is equivalent to the
existence of an equivalent $q$-uniformly smooth norm.

\end{abstract}

 \setcounter{section}{-1}


 \makeatletter
 \renewcommand{\@makefntext}[1]{#1}

 \makeatother \footnotetext{\noindent
  C.O.: Wuhan Institute of Physics and Mathematics, Chinese Academy of Sciences
  P.O. Box 71010\\
  Wuhan 430071,  China\\
  ouyang@wipm.ac.cn\\
  Q.X.: Laboratoire de Math{\'e}matiques, Universit{\'e} de France-Comt{\'e},
  25030 Besan\c{c}on Cedex,  France\\
  qxu@univ-fcomte.fr\\
  2000 {\it Mathematics subject classification:}
  46E40, 42B25,  46B20\\
{\it Key words and phrases}: BMO, Carleson measures, Lusin type,
Lusin cotype, uniformly convex spaces, uniformly smooth spaces.}


\section{Introduction}


Let $\T$ be the unit circle of the complex plane equipped with
normalized Haar measure $dm$. Recall that an integrable function
$f$ on $\T$ is of bounded mean oscillation (BMO) if
 $$\|f\|_{*}=\sup_I\frac1{|I|}\,\int_I|f-f_I|dm<\8,$$
where the supremum runs over all arcs of $\T$ and
 $f_I=|I|^{-1}\int_I fdm$
is the mean of $f$ over $I$. Let $BMO(\T)$ denote the space of BMO
functions on $\T$. The means over arcs in this definition can be
replaced by the averages of $f$  against the Poisson kernel
$P_{z_0}$ for the unit disc $D$:
 $$P_{z_0}(z)=\frac{1-|z_0|^2}{|1-\bar{z_0}z|^2}\,,\quad
 z_0\in D,\, z\in \T.$$
Then
 $$\|f\|_{*}^2\approx
 \sup_{z_0\in D}\int_{\T}|f(z)-f(z_0)|^2P_{z_0}(z)dm(z)$$
with universal equivalence constants. Here as well as in the
sequel, we denote also by $f$ its Poisson integral in $D$:
 $$f(z_0)=\int_{\T} f(z)P_{z_0}(z)dm(z),\quad z_0\in D.$$
On the other hand, it is classical that BMO functions can be
characterized by Carleson measures. A positive measure $\mu$ on
$D$ is called a Carleson measure if
 $$\|\mu\|_C=
 \sup_{z_0\in D}\,\int_{D}\frac{1-|z_0|^2}{|1-\bar{z_0}z|^2}\,
 d\mu(z)<\8.$$
Let $f\in L^1(\T)$. Then $f\in MBO(\T)$ iff
 $|\nabla f(z)|^2\, (1-|z|^2)\,dA(z)$
is a Carleson measure, where $dA(z)$ denotes Lebesgue measure on
$D$. In this case, we have
 \beq\label{bc}
 \|f\|_{*}^2\approx \sup_{z_0\in D}\int_{D}
 |\nabla f(z)|^2\,
 \frac{(1-|z|^2)(1-|z_0|^2)}{|1-\bar{z_0}z|^2}\,dA(z).
 \eeq
We refer to \cite{garnett} for all these results.

This paper concerns the vector-valued version of \eqref{bc}. More
precisely, we are interested in characterizing Banach spaces $B$
for which one of the two inequalities in \eqref{bc} holds for
$B$-valued functions $f$. Given a Banach space $B$ let $L^p(\T;
B)$ denote the usual $L^p$-space of Bochner $p$-integrable
functions on $\T$ with values in $B$. The space $BMO(\T; B)$ of
$B$-valued functions on $\T$ is defined in the same way as in the
scalar case just by replacing the absolute value of $\com$ by the
norm of $B$. Then the vector-valued analogue of \eqref{bc} is the
following:
 \beq\label{bcv}
 c_1^{-1}\,\|f\|_{*}^2\le \sup_{z_0\in D}\int_{D}
 \|\nabla f(z)\|^2\, \frac{(1-|z|^2)(1-|z_0|^2)}{|1-\bar{z_0}z|^2}\,dA(z)
 \le c_2\,\|f\|_{*}^2
 \eeq
for all $f\in BMO(\T; B)$, where $c_1$, $c_2$ are two positive
constants (depending on $B$), and where
 $$\|\nabla f(z)\|=
 \big\|\frac{\partial f}{\partial x}(z)\big\|
 + \big\|\frac{\partial f}{\partial
 y}(z)\big\|\,,\quad z=x+iy.$$
It is part of  the folklore that \eqref{bcv} holds  iff $B$ is
isomorphic to a Hilbert space (see \cite{bla-ill}). We include a
proof of this result at the end of the paper for the convenience
of the reader.

However, if one considers the validity of only one of the two
inequalities in \eqref{bcv}, the matter becomes much subtler and
the corresponding class of Banach spaces is much larger. The
following theorem solves this problem.

\begin{thm}\label{bcvT}
 Let $B$ be a Banach space.
 \begin{enumerate}[\rm (i)]

\item There exists a positive constant $c$ such that
 $$ \sup_{z_0\in D}\int_{D}\|\nabla f(z)\|^2\,
 \frac{(1-|z|^2)(1-|z_0|^2)}{|1-\bar{z_0}z|^2}\,dA(z)
 \le c\,\|f\|_{*}^2$$
holds for all trigonometric polynomials $f$ with coefficients in
$B$ iff $B$ admits an equivalent norm which is $2$-uniformly
convex.

\item There exists a positive constant $c$ such that
 $$ \sup_{z_0\in D}\int_{D}\|\nabla f(z)\|^2\,
 \frac{(1-|z|^2)(1-|z_0|^2)}{|1-\bar{z_0}z|^2}\,dA(z)
 \ge c^{-1}\,\|f\|_{*}^2$$
for all trigonometric polynomials $f$ with coefficients in $B$ iff
$B$ admits an equivalent norm which is $2$-uniformly smooth.

\end{enumerate}
 \end{thm}

We refer to the next section for the definition of uniform
convexity (smoothness). This theorem is intimately related to the
main result of \cite{xu-lp}, where the vector-valued
Littlewood-Paley theory is studied.  Given $f\in L^1(\T; B)$
define the Littlewood-Paley $g$-function
 $$\big(G(f)(z)\big)^2=\int_0^1(1-r)\,\|\nabla
 f(rz)\|^2\, dr\,,\quad z\in\T.$$
The following fact is again well-known: the equivalence
 $$\big\|G(f)\big\|_{L^2(\T)}\approx \big\|f-f(0)\|_{L^2(\T; B)}$$
holds uniformly for all $B$-valued trigonometric polynomials $f$
iff $B$ is isomorphic to a Hilbert space. However, the two
one-sided inequalities are related to uniform convexity
(smoothness). More precisely, we have the following result from
\cite{xu-lp}.

\begin{thm}\label{lusin2}
 Let $B$ be a Banach space.
 \begin{enumerate}[\rm (i)]

\item $B$ has an  equivalent $2$-uniformly convex norm iff for
some $p\in(1,\8)$ $($or equivalently, for every $p\in(1,\8))$
there exists a positive constant $c$ such that
 \beq\label{c2}
 \big\|G(f)\big\|_{L^p(\T)}\le c\,\big\|f\|_{L^p(\T; B)}
 \eeq
holds for all $B$-valued trigonometric polynomials $f$.

\item $B$ has an  equivalent $2$-uniformly smooth norm iff for
some $p\in(1,\8)$ $($or equivalently, for every $p\in(1,\8))$
there exists a positive constant $c$ such that
 \beq\label{t2}
 \big\|f-f(0)\|_{L^p(\T; B)}\le c\,\big\|G(f)\big\|_{L^p(\T)}
 \eeq
holds for all $B$-valued trigonometric polynomials $f$.

\end{enumerate}
 \end{thm}

According to \cite{xu-lp}, the spaces satisfying \eqref{c2} (resp.
\eqref{t2}) are said to be of Lusin cotype 2 (resp. Lusin type 2).
The name Lusin refers to the fact that the Littlewood-Paley
$g$-function can be replaced by the Lusin area function. At this
stage, let us also recall that by Pisier's renorming theorem
\cite{pis-mart}, $B$ has an equivalent $2$-uniformly convex (resp.
smooth) norm iff $B$ is of martingale cotype (resp. type) 2

The value $p=\8$ is, of course, not allowed in Theorem
\ref{lusin2}. At the time of the writing of \cite{xu-lp}, the
second named author guessed that a right substitute of Theorem
\ref{lusin2} for $p=\8$ should be Theorem \ref{bcvT} but could not
confirm this. Our proof of Theorem \ref{bcvT} heavily relies on
Theorem \ref{lusin2} and Calder\'on-Zygmund singular integral
theory. In fact, we will work in the more general setting of an
Euclidean space $\real^n$ instead of $\T$. On the other hand, the
power $2$ in $\|\nabla f\|^2$ plays no longer any special role in
the vector-valued setting. We will consider the analogue of
Theorem \ref{bcvT} for $\|\nabla f\|^q$ with $1<q<\8$. The
corresponding result is stated separately in Theorems \ref{bcc}
and \ref{cbs} below, which correspond to the end point $p=\8$ of
the results of \cite{mtx-adv} and \cite{xu-lp}.


\section{Preliminaries}


Our references for harmonic analysis are \cite{gar-rubio},
\cite{garnett} and \cite{st-ha}. All results quoted in the sequel
without explicit reference can be found there. However, one needs
sometimes to adapt arguments in the scalar case to the
vector-valued setting.

Let $(\O,\mu)$ be a measure space and $B$ a Banach space. For
$1\le p\le\8$ we denote by $L^p(\O,\mu; B)$ the usual $L^p$-space
of Bochner (or strongly) measurable functions on $\O$ with values
in $B$. The norm of $L^p(\O,\mu; B)$ is denoted by $\|\,\|_p$. The
$n$-dimensional Euclidean space $\real^n$ is equipped with
Lebesgue measure. $L^1_{\rm loc}(\real^n;B)$ denotes the space of
locally integrable functions on $\real^n$ with values in $B$.
Recall the Poisson kernel on $\real^n$:
 $$P_t(x)=\frac{\Gamma(\frac{n+1}2)}{\pi^{(n+1)/2}}\,
 \frac{t}{(t^2+|x|^2)^{(n+1)/2}}\,,\quad x\in\real^n\,,\; t>0.$$
Let $f\in L^1_{\rm loc}(\real^n;B)$ such that
 $$\int_{\real^n}\|f(x)\|\,\frac1{1+|x|^{n+1}}\,dx<\8.$$
The Poison integral of $f$ is then defined by
 $$P_t*f(x)=\int_{\real^n}P_t(x-y)f(y)\,dy.$$
The function $P_t*f(x)$ is harmonic in the upper half space
$\real^{n+1}_+$. Let us make a convention  used throughout this
paper: For a function $f\in L^1_{\rm loc}(\real^n;B)$ we also
denote by $f$ its Poisson integral (whenever the latter exists);
thus $f(x, t)=P_t*f(x)$.

The space $BMO(\real^n; B)$ is defined as the space of all
functions $f\in L^1_{\rm loc}(\real^n;B)$ such that
 $$\|f\|_*=\sup_Q\frac1{|Q|}\,\int_Q\|f(x)-f_Q\|\,dx<\8,$$
where the supremum  runs over all cubes $Q\subset\real^n$ (with
sides parallel to the axes), and where $f_Q$ denotes the mean of
$f$ over $Q$. Equipped with $\|\,\|_*$, $BMO(\real^n; B)$ is a
Banach space modulo constants. $BMO(\real^n; \com)$ is simply
denoted by $BMO(\real^n)$.

We will also need the Hardy space $H^1$. There exist several
different (equivalent) ways to define this. It is more convenient
for us to use atomic decomposition. A $B$-valued atom is a
function $a\in L^\8(\real^n; B)$ such that
 $${\rm supp}(a)\subset Q,\quad \int_{\real^n}a\,dx=0,\quad \|a\|_\8\le
 \frac{1}{|Q|}$$
for some cube $Q\subset\real^n$.  We then define $H^1_a(\real^n;
B)$ to be the space of all functions $f$ which can be written as
 $$f=\sum_{k\ge1}\l_k\,a_k$$
with $a_k$ atoms and  $\l_k$ scalars such that $\sum_k|\l_k|<\8$.
The norm of $H^1_a(\real^n; B)$ is defined by
 $$\big\|f\big\|_{H^1_a(\real^n; B)}=\inf\big\{
 \sum_{k\ge1}|\l_k|\,:\, f=\sum_{k\ge1}\l_k\,a_k\big\}.$$
This is a Banach space. It is well-known that $H^1_a(\real^n; B)$
coincides with the space of all $f\in L^1(\real^n; B)$ such that
 $$\sup_{t>0}\|f(\,\cdot, t)\|\in L^1(\real^n).$$
Fefferman's duality theorem between $H^1$ and $BMO$ remains valid
in this setting (with a slight condition on $B$). More precisely,
$BMO(\real^n; B^*)$ is isomorphically identified as a subspace of
the dual $H^1_a(\real^n; B)^*$; moreover, it is norming in the
following sense: for any $f\in H^1_a(\real^n; B)$
 $$\big\|f\big\|_{H^1_a(\real^n; B)}\approx
 \sup\big\{|\la f,\; g\ra|\,:\, g\in BMO(\real^n; B^*),\;
 \|g\|_*\le1\big\}$$
with universal equivalence constants. Note that this duality
result follows immediately from the atomic definition of $H_a^1$.
If $B^*$ has the Radon-Nikodym property (in particular, if $B$ is
reflexive), then
 $$H^1_a(\real^n; B)^*=BMO(\real^n; B^*).$$
We refer to \cite{bla-dua} and \cite{bourg-sin} for more details.

\smallskip

BMO functions can be characterized by Carleson measures. Let
 $$\Ga=\{(z, t)\in\real^{n+1}_+\,:\, |z|<t\},$$
the standard cone of $\real^{n+1}_+$. $\Ga(x)$ denotes the
translation of $\Ga$ by $(x, 0)$ for $x\in\real^n$: $\Ga(x)=\Ga
+(x, 0)$. Let $Q$ be a cube. The tent over $Q$ is defined by
 $$\wh Q=\real^{n+1}_+\setminus\bigcup_{x\in Q^c}\Ga(x)\,.$$
A positive measure $\mu$ on $\real^{n+1}_+$ is called a Carleson
measure if
 $$\|\mu\|_C=\sup_{Q\; {\rm cube}}\,\frac{\mu(\wh Q)}{|Q|}<\8.$$
Then $f\in BMO(\real^n)$ iff $\mu(f)=(t|\nabla f(x, t)|)^2dxdt/t$
is a Carleson measure. Moreover,
 $$
 \|f\|_*^2\approx \|\mu(f)\|_C.
 $$
This is the analogue of \eqref{bc} for $\real^n$. Our main concern
is the validity of each of the two one-sided inequalities of the
equivalence above in the vector-valued setting.  The previous
result is, of course, part of the Littlewood-Paley theory. In this
regard let us recall its $L^p$-analogue. Let $f\in L^p(\real^n)$.
Define the Lusin integral function of $f$:
 $$\big(S(f)(x)\big)^2=
 \int_{\Ga}(t|\nabla f(x+z, t)|)^2\,\frac{dzdt}{t^{n+1}}\,,\quad
 x\in\real^n\,.$$
Then
 \be
 \|f\|_p\approx \|S(f)\|_p\,,\quad
 \forall\; f\in L^p(\real^n),\; 1<p<\8.
 \ee
The vector-valued Littlewood-Paley theory is studied in
\cite{mtx-adv} and \cite{xu-lp}. Let $1<q<\8$ and $f\in
L^p(\real^n; B)$. Define
 $$\big(S_q(f)(x)\big)^q=
 \int_{\Ga}(t\|\nabla f(x+z, t)\|)^q\,\frac{dzdt}{t^{n+1}}\,,\quad
 x\in\real^n\,,$$
where
 $$\|\nabla f(x, t)\|=
 \big\|\frac{\partial}{\partial t}f(x,t)\big\| +
 \sum_{i=1}^n \big\|\frac{\partial}{\partial x_i}f(x,t)\big\|.$$
According to \cite{xu-lp} and \cite{mtx-adv}, $B$ is said to be of
Lusin cotype $q$ if for some $p\in(1,\8)$ (or equivalently, for
every $p\in(1,\8)$) there exists a positive constant $c$ such that
 $$\big\|S_q(f)\big\|_{p}\le c\,\big\|f\|_{p}$$
for all compactly supported $B$-valued continuous  functions $f$
on $\real^n$. Similarly, we define Lusin type $q$ by reversing the
inequality above. Note that if $B$ is of Lusin cotype (resp. type)
$q$, then necessarily $q\ge 2$ (resp. $q\le2$). By \cite{xu-lp}
and \cite{mtx-adv}, Lusin cotype (resp. type) $q$ is equivalent to
martingale cotype (resp. type) $q$. We will not need the latter
notion and refer the interested reader to \cite{pis-proba} and
\cite{pis-mart}.   By Pisier's renorming theorem \cite{pis-mart},
$B$ is of martingale cotype (resp. type) $q$ iff $B$ has an
equivalent norm which is $q$-uniformly convex (resp. smooth). Let
us recall this last notion for which we refer to \cite{LT-II} for
more information. First define the modulus of convexity and
modulus of smoothness of $B$ by
 \be
 \d_B(\e)
 &=&\inf\big\{1-\big\|\frac{a+b}2\big\|\;:\; a,\; b\in B,\;
 \|a\|=\|b\|=1,\; \|a-b\|=\e\big\},\quad 0<\e<2,\\
 \rho_B(t)
 &=&\sup\big\{\frac{\|a+tb\|+\|a-tb\|}2\,-1\;:\; a,\; b\in
B,\;
 \|a\|=\|b\|=1\big\},\quad t>0.
 \ee
$B$ is called uniformly convex if $\d_B(\e)>0$ for every $\e>0$,
and uniformly smooth if $\lim_{t\to0}\rho_B(t)/t=0$. On the other
hand, if $\d_B(\e)\ge c\,\e^q$ for some positive constants $c$ and
$q$, $B$ is called $q$-uniformly convex. Similarly, we define
$q$-uniformly smoothness by demanding $\rho_B(t)\le c\, t^q$ for
some $c>0$ and $q>1$. It is well-known that for $1<p<\8$ any
(commutative or noncommutative) $L^p$-space is $\max(2,
p)$-uniformly convex and $\min(2, p)$-uniformly smooth.


\section{A singular integral operator}
 \label{A singular integral operator}

Let the cone $\Ga=\{(z, t)\in\real^n_+\,:\, |z|<t\}$ be equipped
with the measure $dzdt/t^{n+1}$. Let $1<q<\8$ and $B$ be a Banach
space. Set $A=L^q(\Ga; B).$ For $h\in L^p(\real^n; A)$ we will
consider $h$ as a function of either a sole variable $x\in\real^n$
or three variables $(x, z,t)\in\real^n\times\Ga$. In the first
case $h(x)$ is a function of two variables $(z,t)$ for every
$x\in\real^n$. Thus $h(x)(z,t)=h(x, z,t)$.

We will consider singular integral operators with kernels taking
values in $\L(A)$, the space of bounded linear operators on $A$.
Recall that $P_t$ denotes the Poisson kernel on $\real^n$. Let
 \beq\label{phi}
 \f_t(x)=t\frac{\partial}{\partial t}P_t(x).
 \eeq
For $h\in L^p(\real^n; A)$  define
 \beq\label{Phi}
 \Phi(h)(x, u, s)=\int_{\Ga}\int_{\real^n}
 \f_s*\f_t(x+u+z-y)h(y, z,t)\,dy\,\frac{dzdt}{t^{n+1}}.
 \eeq
$\Phi(h)$ is well defined for $h$ in a dense vector subspace of
$L^p(\real^n; A)$. Indeed, let $h: \real^n\to A$ be a compactly
supported continuous function  such that for each $x\in {\rm
supp}(h)$ the function $h(x): \Ga\to B$ is continuous and
supported by a compact of $\Ga$ independent of $x$. Then it is
easy to check that $\Phi(h)$ is well defined and belongs to
$L^p(\real^n;A)$ for all $p$. On the other hand, it is clear that
the family of all such functions $h$ is dense in $L^p(\real; A)$
for every $p<\8$. In the sequel, $h$ will be assumed to belong to
this family whenever we consider $\Phi(h)$.

The following will be crucial later on. We refer to  \cite{xu-lp}
for a similar lemma on the circle $\T$.

\begin{lem}\label{sing}
 The map $\Phi$ extends to a bounded map on $L^p(\real^n; A)$
for every $1<p<\8$, and also a bounded map from $H^1_a(\real^n;
A)$ to $L^1(\real^n; A)$. Moreover, denoting again by $\Phi$ the
extended maps, we have
 $$\|\Phi: L^p(\real^n; A)\to L^p(\real^n; A)\|\le c\,,\quad
 \|\Phi: H^1_a(\real^n; A)\to L^1(\real^n; A)\|\le c\,,$$
where the constant $c$ depends only on $p, q$ and $n$.

A similar statement holds for each of the $n$ partial derivatives
in $x_i$  instead  of $\partial/\partial t$ in the definition of
$\f$ in \eqref{phi}.
 \end{lem}

\pf The proof is based on  Calder\'on-Zygmund singular integral
theory for vector-valued kernels for which we refer to
\cite{gar-rubio}. We will represent $\Phi$ as a singular integral
operator. Let
 $$k_{s,t}(x)=\f_s*\f_t(x)=\int_{\real^n}\f_s(x-y)\f_t(y)dy.$$
Then
 \beq\label{Phik}
 \Phi(h)(x, u, s)=\int_{\Ga}\int_{\real^n}
 k_{s,t}(x+u+z-y)h(y, z,t)\,dy\,\frac{dzdt}{t^{n+1}}.
 \eeq
On the other hand, using the definition of $\f_t$ and the
semigroup property of $P_t$, we find
 \beq\label{kerk}
 k_{s,t}(x)=st\frac{\partial^2}{\partial
 r^2}P_r(x)\big|_{r=s+t}\,.
 \eeq
Now consider the operator-valued kernel $K(x): A\to A$ defined by
 \be
  K(x)(a)(u,s)=\int_{\Ga}k_{s,t}(x+u+z)a(z,
  t)\frac{dzdt}{t^{n+1}}\,,\quad a\in A.
  \ee
Then $\Phi(h)$ can be rewritten as
 $$\Phi(h)(x)=K*h(x)=\int_{\real^n}K(x-y)(h(y))dy.$$
Thus $\Phi$ is a convolution operator with kernel $K$. We will
show that $K$ is a regular Calder\'on-Zygmund kernel with values
in $\L(A)$. Namely, $K$ satisfies the following norm estimates
 $$\|K(x)\|\le \frac{c}{|x|^{n}}\quad\mbox{and}\quad
 \|\nabla K(x)\|\le \frac{c}{|x|^{n+1}}$$
for some positive constant $c$ depending only $n$. To this end
first observe that by \eqref{kerk}
 \beq\label{kestimate}
 |k_{s, t}(x)|\le \frac{cst}{(s+t+|x|)^{n+2}}\,.
 \eeq
Here as well as in the rest of the paper, letters $c, c', c_1... $
denote positive constants which may depend on $n$, $q$, $p$ or $B$
but never on particular functions in consideration. They may also
vary from line to line. Let $a\in A$ with $\|a\|\le 1$. Let $q'$
denote the conjugate index of $q$. Then by the H\"older inequality
and \eqref{kestimate}, we deduce
 \be
 \|K(x)(a)(u, s)\|^{q'}\le c^{q'}\,\int_{\Ga}
 \frac{s^{q'}t^{q'}}{(s+t+|x+u+z|)^{(n+2)q'}}\,
 \frac{dzdt}{t^{n+1}}\,.
 \ee
Since $|z|<t$, we have
 $$\frac12\,(s+t+|x+u|)\le s+t+|x+u+z|\le 2(s+t+|x+u|).$$
It then follows that
 \be
 \|K(x)(a)(u, s)\|^{q'}
 &\le& c_1^{q'}\,\int_{\Ga}
 \frac{s^{q'}t^{q'}}{(s+t+|x+u|)^{(n+2)q'}}\,
 \frac{dzdt}{t^{n+1}}\\
 &\le& c_2^{q'}\,\frac{s^{q'}}{(s+|x+u|)^{(n+1)q'}}\,.
 \ee
Therefore,
 \be
 \|K(x)(a)\|_A^{q}
 &=&\int_{\Ga} \|K(x)(a)(u, s)\|^{q}\,\frac{duds}{s^{n+1}}\\
 &\le& c_2^{q}\,\int_{\Ga}\frac{s^{q}}{(s+|x+u|)^{(n+1)q}}
 \,\frac{duds}{s^{n+1}}
 \le \frac{c_3^{q}}{|x|^{nq}}\,.
 \ee
Taking the supremum over all $a$ in the unit ball of $A$, we
deduce that $K(x)$ is a bounded operator on $A$ and
 $$\|K(x)\|\le \frac{c_3}{|x|^{n}}\,.$$
Similarly, we show
 $$\|\nabla K(x)\|\le \frac{c_4}{|x|^{n+1}}\,.$$
Therefore, $K$ is a regular vector-valued kernel.

Since $\Phi$ is the singular integral operator with kernel $K$, by
\cite[Theorem~V.3.4]{gar-rubio} (see also
\cite[Theorem~4.1]{mtx-adv}), the lemma is reduced to the
boundedness of $\Phi$ on $L^p(\real^n; A)$ for some $p\in(1,\8)$.
Clearly, the most convenient choice of $p$ is $p=q$. By
\eqref{Phik} and the H\"older inequality
 $$
 \|\Phi(h)(x, u, s)\|
 \le \a\,\cdot\, \b,$$
where
 \be
 \a^{q'}
 &=& \int_{\Ga}\int_{\real^n}
 |k_{s,t}(x+u+z-y)|\,dy\,\frac{dzdt}{t^{n+1}}\,,\\
 \b^{q}
 &=& \int_{\Ga}\int_{\real^n}
 |k_{s,t}(x+u+z-y)|\,\|h(y, z,t)\|^q\,dy\,\frac{dzdt}{t^{n+1}}\,.
 \ee
Using \eqref{kestimate}, we find
 \be
 \a^{q'}
 &\le& c\, \int_{\Ga}\int_{\real^n}
 \frac{st}{(s+t+|x+u+z-y|)^{n+2}}\,\,dy\,\frac{dzdt}{t^{n+1}}\\
 &\le& c_1\, \int_{\Ga}
 \frac{st}{(s+t)^{2}}\,\,\frac{dzdt}{t^{n+1}}\le c_2\,.
 \ee
Hence,
 \be
 \big\|\Phi(h)\big\|_{q}^q
 &=&\int_{\real^n}\int_{\Ga} \|\Phi(h)(x, u, s)\|^q\,
 \frac{duds}{s^{n+1}}\,dx\\
 &\le& c_3\int_{\real^n}\int_{\Ga}\int_{\Ga}\int_{\real^n}
 \frac{st}{(s+t+|x+u+z-y|)^{n+2}}\,dx\,\frac{duds}{s^{n+1}}\,
 \|h(y, z, t)\|^q\,\frac{dzdt}{t^{n+1}}\,dy\\
 &\le& c_4 \int_{\real^n}\int_{\Ga}
 \|h(y, z, t)\|^q\,\frac{dzdt}{t^{n+1}}\,dy
 =c_4\,\big\|h\big\|^q_{L^q(\real^n;A)}\,.
 \ee
Thus $\Phi$ extends to a bounded map on $L^q(\real^n;A)$, so the
lemma is proved.\cqd


\section{Carleson measures and uniform convexity}


The following theorem is the main result of this section. Recall
that $\wh Q$ denotes the tent over $Q$ for a cube
$Q\subset\real^n$.

\begin{thm}\label{bcc}
 Let $B$ be a Banach space and $2\le q<\8$. Then the following
statements are equivalent:
 \begin{enumerate}[\rm(i)]

 \item There exists a positive constant $c$ such that
 \beq\label{Cq}
 \left(\sup_{Q\;{\rm cube}}\frac1{|Q|}\,\int_{\wh Q}
 \big(t\,\|\nabla f(x,t)\|\big)^q\, \frac{dxdt}t\right)^{1/q}
 \le c\,\|f\|_*\,,\quad\forall\; f\in BMO(\real^n;B).
 \eeq

 \item $B$ has an equivalent norm which is $q$-uniformly convex.

 \end{enumerate}
 \end{thm}

Inequality \eqref{Cq} means that $\big(t\,\|\nabla
f(x,t)\|\big)^q\, dxdt/t$ is a Carleson measure on $\real^{n+1}_+$
for every $f\in BMO(\real^n; B)$. In this regard, let us introduce
a more function $C_q$, besides the Lusin function $S_q$. Given $f:
\real^n\to B$ define
 \beq\label{def of Cq}
 C_q(f)(x)=\left(\sup_{Q}\frac1{|Q|}\,\int_{\wh Q}
 \big(t\,\|\nabla f(y,t)\|\big)^q\, \frac{dydt}t\right)^{1/q}\,,
 \eeq
where the supremum runs over all cubes $Q$ containing $x$. Then
\eqref{Cq} can be rephrased as
 $$\|C_q(f)\|_\8\le c\,\|f\|_*\,.$$

The proof of Theorem \ref{bcc} and that of Theorem \ref{cbs} below
heavily rely on the results on Lusin type and cotype in
\cite{mtx-adv}. We collect them in the following lemma for the
convenience of the reader and also for later reference.

\begin{lem}\label{mtx}
 Let $B$ be a Banach space and $2\le q<\8$. Then the following
statements are equivalent:
 \begin{enumerate}[\rm(i)]
 \item $B$ is of Lusin cotype $q$. Namely,
for some $p\in(1,\8)$ $($or equivalently, for every $p\in(1,\8))$
there exists a positive constant $c$ such that
 $$\big\|S_q(f)\big\|_{p}\le c\,\big\|f\|_{p}\,,\quad
 \forall\; f\in L^p(\real^n; B).$$
 \item There exists a constant $c$ such that
 $$\big\|S_q(f)\big\|_1\le c\,\big\|f\|_{H^1_a(\real^n; B)}\,,\quad
 \forall\; f\in H^1_a(\real^n; B).$$
 \item  $B$ has an  equivalent $q$-uniformly convex norm.
 \item $B^*$ is of Lusin type $q'$, where $q'$ is the conjugate
 index of $q$.
 \item $B^*$ has an  equivalent $q'$-uniformly smooth norm.
 \end{enumerate}
 \end{lem}

\n{\it Proof of Theorem \ref{bcc}.} (ii)$\;\Rightarrow\;$(i). Let
$f\in BMO(\real^n;B)$ with $\|f\|_*\le 1$. Let $Q\subset \real^n$
be a cube. Set $\wt Q=2Q$, the cube of the same center as $Q$ and
of double side length. Write
 $$f=(f-f_{\wt Q})\un_{\wt Q} + (f-f_{\wt Q})\un_{\wt Q^c} +
 f_{\wt Q}\;
 {\mathop =^{\rm def}}\; f_1+f_2 + f_{\wt Q}\,.$$
Then
 $$\nabla f(x,t)=\nabla f_1(x,t)+\nabla f_2(x,t),$$
so
 $$\left(\frac1{|Q|}\,\int_{\wh Q}
 \big(t\,\|\nabla f(x,t)\|\big)^q\, \frac{dxdt}t\right)^{1/q}
 \le \a_1 + \a_2\,,$$
where
 $$\a_k=\left(\frac1{|Q|}\,\int_{\wh Q}
 \big(t\,\|\nabla f_k(x,t)\|\big)^q\, \frac{dxdt}t\right)^{1/q}\,,
 \quad k=1,\, 2.$$
For $\a_1$ by the Fubini theorem we have
 \be
 |Q|\,\a_1^q
 &\le& c_n^q\int_Q\int_{\Ga} \big(t\,\|\nabla
 f_1(x+z,t)\|\big)^q\,\frac{dzdt}{t^{n+1}}\,dx\\
 &=&c_n^q\int_Q\big(S_q(f_1)(x)\big)^qdx
 \le c_n^q\|S_q(f_1)\|_q^q\,,
 \ee
where $c_n$ is a constant depending only on $n$. By (ii) and Lemma
\ref{mtx}, $B$ is of Lusin cotype $q$. Thus
 $$\|S_q(f_1)\|_q\le c\, \|f_1\|_q\,.$$
However, by the John-Nirenberg theorem
 $$\|f_1\|_q\le c'\,|Q|^{1/q}\,\|f\|_*\le c'\,|Q|^{1/q}\,.$$
It then follows that
 $$\a_1\le c_nc\,c'\,.$$
To deal with $\a_2$ we write
 $$\nabla f_2(x,t)=\int_{\real^n}\nabla P_t(x-y)f_2(y)dy
 =\int_{\wt Q^c}\nabla P_t(x-y)f_2(y)dy.$$
Note that
 $$|\nabla P_t(x-y)|\le \frac{c_n}{(t+|x-y|)^{n+1}}\,.$$
On the other hand, for $(x, t)\in \wh Q$ and $y\in\wt Q^c$
 $$\frac1{(t+|x-y|)^{n+1}}\approx \frac1{(\el+|x-y|)^{n+1}}\,,$$
where $\el=\el(Q)$ is the  side length of $Q$. Thus
 \be
 \|\nabla f_2(x,t)\|
 &\le& c_n'
 \int_{\wt Q^c}\|f_2(y)\|\,\frac1{(\el+|x-y|)^{n+1}}\,dy\\
 &\le& \frac{c_n''}\el\int_{\real^n}\|f_2(y)\|\,P_\el(x-y)\,dy.
 \ee
We now use a well-known characterization of BMO functions, in
which averages over cubes are replaced by averages against the
Poisson kernel.  Namely, a function $g:\real^n\to B$ belongs to
$BMO(\real^n; B)$ iff
 \be
 \sup_{(x, t)\in\real^n_+}
 \int_{\real^n}\|g(y)-g(x,t)\|\,P_t(x-y)\,dy<\8\,.
 \ee
If this is the case, the supremum above is equivalent to $\|g\|_*$
with relevant constants depending only on $n$.  Then we deduce
 \be
 \|\nabla f_2(x,t)\|
 \le \frac{c}\el.
 \ee
 Therefore,
 $$
 \a_2^q
 \le \frac{c^q}{\el^q\,|Q|}\,\int_{\wh Q}
 t^q\, \frac{dxdt}t \le c'\,.$$
Combining the preceding inequalities, we find that
$\big(t\,\|\nabla f(x, t)\|\big)^q dxdt/t$ is a Carleson measure
on $\real^{n+1}_+$ with constant depending only on $n, q$ and $B$
for every $f\in BMO(\real^n; B)$ with $\|f\|_*\le 1$. This
concludes the proof of (ii)$\;\Rightarrow\;$(i).

(i)$\;\Rightarrow\;$(ii). This proof is harder.  Let $A=L^q(\Ga;
B)$ (recall that the cone $\Ga$ is equipped with the measure
$dzdt/t^{n+1}$). Given a function $f\in L^p(\real^n; B)$ define
 $$\S_q(f)(x, z, t)=t\,\frac{\partial}{\partial t}f(x+z, t),
 \quad x\in\real^n,\; (z, t)\in\Ga.$$
We regard $\S_q(f)$ as a function on $\real^n$ with values in $A$.
Then
 $$\|\S_q(f)(x)\|_A=S^t_q(f)(x)\,,$$
where $S^t_q(f)$ is the Lusin integral function of $f$ but using
only the partial derivative in $t$. Also note that
 $$\S_q(f)(x, z,t)=\f_t*f(x+z),$$
where $\f$ is defined by \eqref{phi}. As in section \ref{A
singular integral operator}, $\S_q$ can be represented as a
singular integral operator with a regular kernel taking values in
the space of bounded linear maps from $B$ into $A$ (see
\cite{st-sing} for the scalar case and \cite{xu-lp} for $\T$). By
\cite{mtx-adv} and\cite{xu-lp}, (ii) is equivalent to the
following inequality
 \beq\label{L-B}
 \|\S_q(f)\|_*\le c\,\|f\|_\8\,,\quad \forall\; f\in L^\8(\real^n;
 B).
 \eeq
Note that this inequality is a finite dimensional property.
Namely, if \eqref{L-B} holds for every finite dimensional subspace
$E$ of $B$  in place of $B$ with constant independent of $E$, then
\eqref{L-B} holds for the whole $B$ too. Thus we can assume $\dim
B<\8$ in the rest of the proof. To prove \eqref{L-B} we will use
duality. We first show that (i) implies
 \beq\label{H-B*}
 \|g\|_{H^1_a(\real^n; B^*)}
 \le c\,\|S^t_{q'}(g)\|_1
 \eeq
for all compactly supported continuous functions $g:\real^n\to
B^*$. To this end let $f\in BMO(\real^n; B)$ with $\|f\|_*\le1$.
Then by Plancherel's theorem
 \beq\label{plan}
 \int_{\real^n}\la f(x),\; g(-x)\ra\, dx
 =4\int_{\real^{n+1}_+} \la t\frac{\partial}{\partial t}f(x, t),\;
 t\frac{\partial}{\partial t}g(-x, t)\ra
 \,\frac{dxdt}t\,.
 \eeq
Note that since $\dim B<\8$,  this equality is reduced to the
scalar case, in which it is well-known and immediately follows
from Plancherel's theorem. Let $C^t_q(f)$ denote the function
defined by \eqref{def of Cq} using only the partial derivative in
$t$ . Then by \eqref{Cq} we find
  \be
 \left|\int_{\real^n}\la f(x),\; g(-x)\ra\, dx\right|
 &\le& 4\int_{\real^{n+1}_+} t\big\|\frac{\partial}{\partial t}f(x,
 t)\big\|\,
 t\big\|\frac{\partial}{\partial t}g(-x, t)\big\|
 \,\frac{dxdt}t\\
  &\le& c'\,\int_{\real^n} C^t_q(f)(x)\,S^t_{q'}(g)(-x)dx\\
 &\le& c\,c'\, \|f\|_*\,\|S^t_{q'}(g)\|_1\,,
 \ee
where we have used Theorem~1 (a) of \cite{co-me-st} for the next
to the last inequality.  Note that the inequality there is proved
only for $q=2$; but the arguments can be easily modified to our
situation. Taking the supremum over all $f$ in the unit ball of
$BMO(\real^n; B)$, we obtain \eqref{H-B*}.

Return back to \eqref{L-B}. We use again duality, this time that
between $BMO(\real^n; A)$ and $H^1_a(\real^n; A^*)$. Fix a
function $f\in L^\8(\real; B)$. Recall that $\S_q(f)$ is a
function from $\real^n$ to $A$ and the left hand side of
\eqref{L-B} is $\|\S_q(f)\|_{BMO(\real^n; A)}$. Thus it suffices
to prove
 \beq\label{L-B1}
 |\la \S_q(f),\; h\ra|\le c\,\|f\|_\8\, \|h\|_{H^1_a(\real^n;
 A^*)}\,,\quad\forall\; h\in H^1_a(\real^n; A^*).
 \eeq
Again by approximation, we need only to consider a nice $h$. We
have
 \be
 \la \S_q(f),\; h\ra
 &=&\int_{\real^n}\int_{\Ga}\la\f_t*f(x+z),\; h(-x,
 z,t)\ra\,\frac{dzdt}{t^{n+1}}\,dx\\
 &=&\int_{\real^n}\int_{\Ga}\la f(y),\; \f_t(\,\cdot
 +z)*h(\,\cdot, z,t)(-y)\ra\,\frac{dzdt}{t^{n+1}}\,dy\\
 &=&\int_{\real^n}\la f(y),\;\Psi(h)(-y)\ra\,dy,
 \ee
where
 \begin{eqnarray}\label{Psi}
  \begin{array}{ccl}
  \begin{displaystyle}\Psi(h)(x)\end{displaystyle}
  &=& \begin{displaystyle}
  \int_{\Ga}\f_t(\,\cdot+z)*h(\,\cdot,
 z,t)(x)\,\frac{dzdt}{t^{n+1}} \end{displaystyle} \\
  &=& \begin{displaystyle}
  \int_{\Ga}\int_{\real^n}
 \f_t(x+z-y)h(y, z,t)\,dy\,\frac{dzdt}{t^{n+1}}\,.\end{displaystyle}
 \end{array}
 \end{eqnarray}
$\Psi(h)$ is a function on $\real^n$ with values in $B^*$.
Therefore, by \eqref{H-B*}
 \be
 |\la \S_q(f),\; h\ra|
 &\le& \|f\|_\8\, \|\Psi(h)\|_1
 \le\|f\|_\8\, \|\Psi(h)\|_{H_a^1(\real^n; B^*)}\\
 &\le& c\,\|f\|_\8\, \|S^t_{q'}(\Psi(h))\|_1
 =c\,\|f\|_\8\, \|\S_{q'}(\Psi(h))\|_{L^1(\real^n; A^*)}\,.
 \ee
Here we use the same notation $\S$ in the dual setting, which is
consistent with the preceding meaning for $A^*$ is the space
associated to $B^*$ in the same way as $A$ associated to $B$:
 $$A^*=L^{q'}(\Ga; B^*).$$
Now it is easy to see that
 $$\S_{q'}(\Psi(h))=\Phi(h),$$
where $\Phi$ is the map defined by \eqref{phi} with $(B^*, q')$
instead of $(B, q)$. Thus by Lemma \ref{sing},
$$\|\S_{q'}(\Psi(h))\|_{L^1(\real^n; A^*)}\le c\,
 \|h\|_{H^1_a(\real^n;
 A^*)}\,.$$
Combining the preceding  inequalities, we obtain \eqref{L-B1}, and
consequently, \eqref{L-B} too. This shows the implication
(i)$\;\Rightarrow\;$(ii). Thus the proof of the theorem is
complete.\cqd

\medskip

The previous proof of (i)$\;\Rightarrow\;$(ii) shows the following
result, which extends \cite[Theorem~5.3]{mtx-adv} (and
\cite[Theorem~2.5]{xu-lp}) to the case $p=1$.

\begin{cor}
 Let $B$ be a Banach space and $1<q\le2$. Then the following statements are equivalent:
 \begin{enumerate}[\rm(i)]
 \item $B$ is of Lusin type $q$.
 \item There exists a constant $c$ such that
 $$\|f\|_{H^1_a(\real^n; B)}\le c\, \|S_q(f)\|_1$$
holds for all compactly supported continuous functions $f$ from
$\real^n$ to $B$.
 \item There exists a constant $c$ such that
 $$\|f\|_1\le c\, \|S_q(f)\|_1$$
holds for all compactly supported continuous functions $f$ from
$\real^n$ to $B$.
\end{enumerate}
\end{cor}


\section{Carleson measures and uniform smoothness}


This section deals with the properties dual to those in Theorem
\ref{bcc}. The following theorem gives the characterization of
Lusin type in terms of Calrseon measures.

\begin{thm}\label{cbs}
 Let $B$ be a Banach space and $1<q\le2$. Then the following
statements are equivalent:
 \begin{enumerate}[\rm(i)]

 \item There exists a positive constant $c$ such that
 \beq\label{qC}
 \|f\|_* \le c\, \left(\sup_{Q\;{\rm cube}}
 \frac1{|Q|}\,\int_{\wh Q} \big(t\,\|\nabla f(x,t)\|\big)^q\,
 \frac{dxdt}t\right)^{1/q}
 \eeq
holds for all compactly supported continuous functions $f$ from
$\real^n$ to $B$.

 \item $B$ has an equivalent $q$-uniformly smooth norm.

 \end{enumerate}
 \end{thm}

\pf (ii)$\;\Rightarrow\;$(i) First note that by  Lemma \ref{mtx},
(ii) is equivalent to
 \beq\label{Sq'}
 \|S_{q'}(g)\|_1\le c\,\|g\|_{H^1_a(\real^n; B^*)}\,,
 \quad\forall\;g\in
 H^1_a(\real^n; B^*).
 \eeq
Let $f: \real^n\to B$ be a compactly supported continuous
function. We are going to prove \eqref{qC}. This proof is similar
to that of \eqref{H-B*} but in a converse direction.  By
approximation, we can assume that $f$ takes values in a finite
dimensional subspace of $B$; then replacing $B$ by this subspace,
we can simply assume $\dim B<\8$. Using the duality between
$BMO(\real^n; B)$ and $H^1_a(\real^n; B^*)$, we find a function
$g\in H^1_a(\real^n; B^*)$ of unit norm such that
 $$\|f\|_*\approx \int_{\real^n}\la f(x),\;
 g(-x)\ra\, dx,$$
where the equivalence constants depend only on $n$. By
approximation, we can further assume that $g$ is sufficiently nice
so that all calculations below are legitimate.  By Plancherel's
theorem
 \be
 \int_{\real^n}\la f(x),\; g(-x)\ra\, dx
 =\int_{\real^{n+1}_+} \la t\nabla f(x,t),\; t\nabla g(-x, t)\ra
 \,\frac{dxdt}t\,.
 \ee
By \cite{co-me-st} and \eqref{Sq'}, we find
 \be
 \int_{\real^n}\la f(x),\; g(-x)\ra\, dx
 &\le& \int_{\real^{n+1}_+} t\|\nabla f(x,t)\|\, t\|\nabla g(-x, t)\|
 \,\frac{dxdt}t\\
 &\le& c'\int_{\real^n} C_q(f)(x)\,S_{q'}(g)(-x)dx\\
 &\le& c'\|C_q(f)\|_\8\, \|S_{q'}(g)\|_1\\
  &\le& c''\|C_q(f)\|_\8\,\|g\|_{H^1_a(\real^n; B^*)}
  \le c''\|C_q(f)\|_\8\,.
 \ee
Combining the preceding inequalities, we deduce \eqref{qC}.

(i)$\;\Rightarrow\;$(ii). Assume (i). It suffices to  prove
\eqref{Sq'}. We will do this only for the Lusin function involving
the partial derivative $\partial/\partial t$, the others being
treated similarly. Thus let $S_{q'}^t$ denote this Lusin function.
Our task is to show
 \beq\label{Sq'1}
  \|S^t_{q'}(g)\|_1\le c\,\|g\|_{H^1_a(\real^n; B^*)}\,,\quad
  \forall\; g\in H^1_a(\real^n; B^*).
  \eeq
We can clearly assume  $\dim B<\8$. Let $A=L^q(\Ga; B)$ be as in
section \ref{A singular integral operator} and keep the notations
introduced there. Note that
 $$A^*=L^{q'}(\Ga; B^*).$$
Now fix a nice function $g\in H^1_a(\real^n; B^*)$. Recall that
 \be \|S^t_{q'}(g)\|_1
 =\int_{\real^n}\left(\int_{\Ga}\big(t\|
 \frac{\partial}{\partial t}g(x+z,t)\|_{B^*}\big)^{q'}\,
 \frac{dzdt}t\right)^{1/q'}\,dx
 =\|\wt g\|_{L^1(\real^n; A^*)}\,,
 \ee
where
 $$\wt g(x, z, t)=t\frac{\partial}{\partial t}
 g(x+z,t).$$
Thus there exists a function $h\in L^\8(\real^n; A)$ of norm $1$
such that
 \be
 \|S^t_{q'}(g)\|_1
 &=& \int_{\real^{n}}\int_{\Ga}\la t\frac{\partial}{\partial t}
 g(x+z,t),\; h(-x,z, t)\ra\,\frac{dzdt}t\,dx\\
 &=&\int_{\real^n}\la g(x),\; \Psi(h)(-x)\ra\,dx,
 \ee
where $\Psi$ is defined by \eqref{Psi}. Therefore, by \eqref{qC},
we deduce
 $$\|S^t_{q'}(g)\|_1
 \le c_n\, \|g\|_{H^1_a(\real^n; B^*)}\, \|\Psi(h)\|_*
 \le c_nc\,\|g\|_{H^1_a(\real^n; B^*)}\,\|C_q(\Psi(h))\|_\8\,.$$
Thus we are reduced to proving
 \be
 \|C_q(\Psi(h))\|_\8\le c\,\|h\|_{L^\8(\real^n; A)}\,,
 \quad \forall\; h\in L^\8(\real^n; A).
 \ee
We will do this only for the partial derivative in the time
variable in the gradient. Namely, we have to show
 \beq\label{C-Ca}
 \frac1{|Q|}\,\int_{\wh Q} \big(s\,\|\frac{\partial}{\partial s}\,
 \Psi(h)(x,s)\|\big)^q\,
 \frac{dxds}s\le c^q\|h\|^q_{L^\8(\real^n; A)}
 \eeq
for any cube $Q\subset\real^n$. The argument below is similar to
the proof of (ii)$\,\Rightarrow\;$(i) in Theorem \ref{bcc}. Using
$\f$ and $k_{s,t}$ in section~\ref{A singular integral operator},
we have
 \be
 s\,\frac{\partial}{\partial s}\,\Psi(h)(x,s)
 =\int_{\real^n}\int_{\Ga}k_{s, t}(x+z-y)h(y, z, t)\,
 \frac{dzdt}t\,dy
 \;{\mathop =^{\rm def}}\;\wt \Phi(h)(x, s).
 \ee
Now fix a cube $Q$ and  a nice $h\in L^\8(\real^n; A)$ with
$\|h\|_{L^\8(\real^n; A)}\le1$. Let $\wt Q=2Q$. Decompose $h$:
 $$h=h\un_{\wt Q} + h\un_{\wt Q^c}
 \;{\mathop =^{\rm def}}\;h_1+h_2.$$
Then \eqref{C-Ca} is reduced to
 $$\b_k=\left(\frac1{|Q|}\,
 \int_{\wh Q} \big(\|\wt\Phi(h_k)(x, s)\|\big)^q\,
 \frac{dxds}s\right)^{1/q}\le c,\quad k=1,\;2.$$
$\b_1$ is easy to estimate. Indeed, using the map $\Phi$ in
\eqref{Phi} and Lemma \ref{sing}, we find
 \be
 |Q|\,\b_1^q
 &\le& c_n^q\int_{Q}\|\Phi(h_1)(x)\|_{A}^q\,dx
 \le c_n^q\|\Phi(h_1)\|_{L^q(\real^n; A)}^q\\
 &\le& c_n^qc^q\|h_1\|_{L^q(\real^n; A)}^q
 \le c_n^qc^q|Q|;
 \ee
whence the desired result for $\b_1$. For $\b_2$ a little more
effort is needed. By \eqref{kestimate}, we have
 \be
 \|\wt\Phi(h_2)(x,s)\|
 \le c\,\int_{\wt Q^c}\int_{\Ga}
 \frac{st}{(s+t+|x+z-y|)^{n+2}}\,\|h(y,z,t)\|\,
 \frac{dzdt}{t^{n+1}}\,dy.
 \ee
By the H\"older inequality and the assumption that
$\|h\|_{L^\8(\Ga; A)}\le1$, the internal integral is estimated as
follows:
 \be
 &&\int_{\Ga}
 \frac{st}{(s+t+|x+z-y|)^{n+2}}\,\|h(y,z,t)\|\,
 \frac{dzdt}{t^{n+1}}\\
 &&~~~\le \left(\int_{\Ga}
 \frac{(st)^{q'}}{(s+t+|x+z-y|)^{(n+2)q'}}\, \frac{dzdt}{t^{n+1}}
 \right)^{1/q'}\,\|h(y)\|_A\\
 &&~~~\le \left(\int_{\Ga}
 \frac{s^{q'}t^{q'}}{(s+t+|x+z-y|)^{(n+2)q'}}\, \frac{dzdt}{t^{n+1}}
 \right)^{1/q'}\\
 &&~~~\approx \frac{s}{(s+|x-y|)^{n+1}}\,.
 \ee
On the other hand, for $(x, s)\in\wh Q$ and $y\in\wt Q^c$, we have
 $$\frac{s}{(s+|x-y|)^{n+1}}\approx \frac{s}{|x-y|^{n+1}}\,.$$
Therefore,
 \be
 \|\wt\Phi(h_2)(x,s)\|
 \le c's\,\int_{\wt Q^c}\frac{dy}{|x-y|^{n+1}}
 \le\frac{c''s}\el,
 \ee
where $\el$ is the side length of $Q$. It then follows that
 $\b_2\le c.$
Thus  \eqref{C-Ca} is proved. This finishes the proof of
\eqref{Sq'1}, so the implication (i)$\;\Rightarrow\;$(ii) too.
\cqd

\medskip

\n{\it Proof of Theorem \ref{bcvT}.} Except the difference between
$\T$ and $\real$, Theorem \ref{bcvT} is a particular case of
Theorems \ref{bcc} and \ref{cbs}. The proofs of these two latter
theorems can be easily adapted to the case of the circle.\cqd

\begin{rk}
 The two ``if'' parts in Theorem \ref{bcvT} can be also proved by
using the invariance of the expression $\|\nabla f(z)\|^2\,dA(z)$
under M\"obius transformations of $D$. This invariance means that
if $w=\f(z)$ is a M\"obius transformation of $D$, then
 $$\|\nabla f(\f(z))\|^2dA(z)=\|\nabla f(w)\|^2dA(w).$$
Now assume that $B$ is $2$-uniformly convex. Then $B$ is of Lusin
cotype $2$. Therefore there exists a constant $c$ such that
 $$\int_{\T}\int_0^1(1-r)\|\nabla f(rz)\|^2dr\,dm(z)
 \le c\, \|f-f(0)\|_2^2\,,\quad\forall\; f\in L^2(\T; B).$$
Then one easily deduces that (with a different $c$)
 $$\int_{D}(1-|z|^2)\|\nabla f(z)\|^2\,dA(z)
 \le c\, \|f-f(0)\|_2^2\,.$$
Now let $z_0\in D$ and let
 $$\f(z)=\frac{z+z_0}{1+\bar{z_0}z}\,.$$
Applying the preceding inequality to $f\circ\f$, we get
 $$\int_{D}\|\nabla f\circ\f(z)\|^2(1-|z|^2)\,dA(z)
 \le c\, \|f\circ\f-f\circ\f(0)\|_2^2\,.$$
Then a change of variables and the previous M\"obius invariance
yield
 $$\int_{D}\|\nabla f(z)\|^2\,
 \frac{(1-|z|^2)(1-|z_0|^2)}{|1-\bar{z_0}z|^2}\,dA(z)
 \le c\,\int_{\T}\|f(z)-f(z_0)\|^2
 P_{z_0}(z)dm(z).$$
Taking the supremum over all $z_0\in D$ gives the first inequality
in Theorem \ref{bcvT}. The same argument applies to the ``if''
part in (ii) there. Unfortunately, this simple proof works neither
for the case of $q\neq2$  nor for that of $\real^n$.

\end{rk}

We end the paper with some comments on \eqref{bcv}. If \eqref{bcv}
holds, then $B$ has an equivalent $2$-uniformly convex norm as
well as an equivalent $2$-uniformly smooth norm. In particular, it
is  of both cotype $2$ and  type $2$, so  isomorphic to a Hilbert
space by Kwapie{\'n}'s theorem \cite{kw} (see also \cite{pis-fact}
to which we refer for the notion of type and cotype too).
Conversely, if $B$ is isomorphic to a Hilbert space, we get
\eqref{bcv} as in the scalar case. Let us give a much more
elementary argument showing that the validity of \eqref{bcv}
implies the isomorphism of $B$ to a Hilbert space. The main point
is the following remark.

\begin{rk}
 Let $1<q<\8$ and $B$ be a Banach space.
Given a finite sequence $(a_k)\subset B$ consider the function
  $$f(z)=\sum_{k\ge1}a_k z^{2^k}\,.$$
Then
 \beq\label{lacu}
 \sup_{z_0\in D}\,\int_{D}(1-|z|)^{q-1}\|f'(z)\|^q\,
 P_{z_0}(z)\,dA(z)\approx
 \sum_{k\ge1}\|a_k\|^q\,
 \eeq
with universal equivalence constants.
 \end{rk}

Recall the following well-known (and easily checked) fact
 $$\|f\|_*\approx
 \big\|\sum_{k\ge1}a_k z^{2^k}\big\|_1\,.$$
Combining this with \eqref{lacu} we deduce the following result
from \cite{bla-ill}: If
 $$\sup_{z_0\in D}\,\int_{D}(1-|z|)^{q-1}\|f'(z)\|^q\,
 P_{z_0}(z)\,dA(z)\le c^q \|f\|_*^q$$
holds for any lacunary polynomial $f$ with coefficients in $B$
with some positive constant $c$, then $B$ is of cotype $q$; the
converse inequality implies that  $B$ is of type $q$.

\medskip

Let us show \eqref{lacu}. Since
 $$f'(z)=\sum_{k\ge1}2^ka_k z^{2^k-1}\,,$$
replacing $a_k$ by $2^ka_k$, we see that \eqref{lacu} is reduced
to
 $$\sup_{z_0\in D}\,\int_{D}(1-|z|)^{q-1}\|f(z)\|^q\,
 P_{z_0}(z)\,dA(z)\approx
 \sum_{k\ge1}2^{-qk}\,||a_k\|^q\,.$$
The lower estimate is very easy. Indeed,  we have (with $z_0=0$)
 \be
 \int_{D}(1-|z|)^{q-1}\|f(z)\|^q\,dA(z)
 &=&\int_0^1(1-r)^{q-1}\int_{\T}\|f(rz)\|^q\,dm(z)\,rdr\\
 &=&\sum_{n\ge1}\int_{1-2^{-n+1}}^{1-2^{-n}}
 (1-r)^{q-1}\int_{\T}\|f(rz)\|^q\,dm(z)\,rdr\\
 &\ge&\sum_{n\ge1}\int_{1-2^{-n+1}}^{1-2^{-n}}
 (1-r)^{q-1}\,\|a_n\|^q\,r^{q2^n}\,rdr\\
 &\approx& \sum_{n\ge1}2^{-qn}\,\|a_n\|^q\,.
 \ee
For the upper estimate, we first majorize $f$ pointwise. For
$n\ge1$ and  $1-2^{-n+1}\le |z|<1-2^{-n}$ we find
 \be
 \|f(z)\|
 \le\sum_{k\le n}\|a_k\| + \sum_{k> n}\|a_k\|\exp(-2^{k-n})\,.
 \ee
Let $0<\a<1$. Then
 $$\sum_{k\le n}\|a_k\|\le c\, 2^{n\a}\,\big(\sum_{k\le n}2^{-k\a
 q}\|a_k\|^q\big)^{1/q}\,.$$
Similarly, for $\b>1$
 \be
 \sum_{k> n}\|a_k\|\exp(-2^{k-n})
 &\le& \big(\sum_{k> n}2^{-k\b q}\|a_k\|^q\big)^{1/q}\,
 \big(\sum_{k>n}2^{k\b q'}\exp(-q'2^{k-n})\big)^{1/q'}\\
 &\le&  c\, 2^{n\b}\,\big(\sum_{k> n}2^{-k\b
 q}\|a_k\|^q\big)^{1/q}\,.
 \ee
It follows that for any $z_0\in D$
 \be
 \int_{D}(1-|z|)^{q-1}\|f(z)\|^q\,
 P_{z_0}(z)\,dA(z)
 &\le& c\,\sum_{n\ge1}2^{-nq}\Big[2^{nq\a}
 \sum_{k\le n}2^{-k\a q}\|a_k\|^q
 +2^{nq\b}\sum_{k> n}2^{-k\b q}\|a_k\|^q\Big]\\
 &\le& c\, \sum_{k\ge1}2^{-qk}\,\|a_k\|^q\,.
 \ee
Therefore, \eqref{lacu} is proved.

\vskip 1cm


\end{document}